\documentclass[11pt]{article}
\usepackage{amsmath}
\usepackage{amssymb}
\usepackage{amsthm}
\usepackage[dvips]{graphicx,color}
\usepackage{psfig}

\newtheorem {theorem}{Theorem}[section]

\newtheorem{lemma}[theorem]{Lemma}

\newtheorem{claim}[theorem]{Claim}

\def\R{{\Bbb R}}

\def\Cov{{\bf{Cov}}}

\def\eps{\epsilon}
\def\bits{{\{0,1\}}}
\def\cube{{\{0,1\}^n}}

\begin{document}

\title{A Law of Large Numbers for Weighted Majority}

\author{Olle H\"{a}ggstr\"{o}m\footnote{Supported by a grant from
the Swedish Research Council.} \\ Chalmers University \\ 
olleh@math.chalmers.se \and 
Gil Kalai\footnote{Supported by an ISF Bikura grant.}\\ 
Hebrew University \\ and Yale University \\
kalai@math.huji.ac.il \and
Elchanan Mossel\footnote{Supported by a Miller fellowship in Computer Science 
and Statistics, U.C. Berkeley.} \\
U.C. Berkeley \\ mossel@stat.berkeley.edu \and}

\date{1 June 2004}
\maketitle

\begin{abstract}

Consider
an election between two candidates  in which  the voters' choices
are random and independent and the probability of a voter choosing
the first candidate is $p>1/2$. Condorcet's Jury Theorem
which he derived from the weak law of large numbers  
asserts that if the number of
voters tends to infinity then  the probability that
the first candidate will be
elected tends to one. 
The notion of influence of a voter or its voting power
is relevant for extensions of the weak law of large numbers 
for
voting rules which are more general than simple majority.
In this paper we point out two different ways to extend 
the classical notions of voting power
and influences to arbitrary probability distributions.
The extension relevant to us is the ``effect'' of a voter, 
which is a weighted version of the correlation 
between the voter's vote and the election's outcomes.  
We prove an extension of the weak law of large numbers 
to weighted majority games 
when all individual effects are small and show that this result 
does not apply to any voting rule which is not based on weighted majority.


\end{abstract}
\section{Introduction}

Consider a biased coin for which the probability for a ``head'' is $p>1/2$.
The weak law of large numbers asserts that if you flip the coin 
$n$ times then the probability that you will see more heads than tails
tends to one as $n$ tends to $\infty$. 
Understanding the scope of the weak law of large numbers 
when the coin flips are 
not independent or when we consider
more complicated events than the event ``to see more heads than tails '',
has attracted considerable attention.

Our motivation came from a game theoretic 
interpretation: Condorcet's
Jury Theorem (see \cite {Y})  asserts that in an election between
two candidates, say  Alice and Bob, if every voter votes for Alice
with probability $p>1/2$ and for Bob with probability $1-p$ and if
these votes are independent, then as the number of voters
tends to infinity the probability that Alice will be elected tends
to one. 
Condorcet's Jury theorem can be interpreted as saying that
even if agents receive very poor (independent) signals indicating
which  decision is correct,  majority voting will nevertheless
result in the correct decision being taken with a  high
probability if  there are enough agents (and each agent votes
according to the  signal he receives). 
This phenomenon is referred to 
as asymptotically complete aggregation of information 
and it plays an important role in theoretical economics.

To describe a more general settings consider the following
framework. Let $f:\cube \to \{0,1\}$ be a Boolean function.
We will assume that $f$ is 
\begin{itemize}
\item
monotone non-decreasing, i.e., 
\[
\left( \forall i: \,\,x_i \geq y_i \right) 
\implies f(x_1,\ldots,x_n) \geq f(y_1,\ldots,y_n), 
\]
\item
and anti-symmetric, i.e., 
$$f(1-x_1,1-x_2,\dots,1-x_n)=1-f(x_1,x_2,\dots,x_n).$$
\end{itemize} 
Let $\mu_p$ denote the product probability measure
on $\cube$ defined by
\begin{equation} \label{eq:product_measure}
\mu_p(x_1, x_2, \dots, x_n) = p^k(1-p)^{n-k},
\end{equation}
where $k=x_1+x_2+\dots+x_n$.
We would like to find 
conditions that guarantee that for a fixed $p>1/2$,  $\mu_p(f)$ is close
to 1. Clearly it is not sufficient that $n$ is large since even
if $f$ is defined
on many variables, it may actually depend only on a few of them.
The notion of {\it influence} of a variable which is 
closely related to notions
of {\it voting power} is important in understanding information aggregation
when we consider general Boolean functions and the product 
probability measure $\mu_p$.
Boolean functions can describe voting rules and are referred to in the 
game theoretic literature as simple games. Anti-symmetric Boolean functions
are called strong simple games.

For a Boolean function $f$ and $x=(x_1,x_2,\dots,x_n) \in \cube$ 
we say that the $k$'th variable
is {\it pivotal} for $f$ if 
$f(x_1,\dots,x_{k-1},0,x_{k+1},\dots,x_n) \ne 
f(x_1,\dots,x_{k-1},1,x_{k+1},\dots,x_n)$.

Let $\mu$ be an arbitrary probability distribution on $\cube$
and let $f$ be a monotone Boolean function that we consider as a voting rule.
Define the 
{\it influence} or the {\it voting power}  of  
of the $k$'th variable 
as the probability
that the $k$'th variable is pivotal. 
Denote by $I_k^{\mu}(f)$ the influence
of the $k$'th variable for the Boolean function $f$, w.r.t. 
the distribution $\mu$.
In other words,
\begin{equation} \label{eq:Idef}
I_k^{\mu}(f) = \mu[(x_1,\ldots,x_n) :  
f(x_1,\ldots,x_{k-1},0,x_{k+1},\ldots,x_n) \neq 
f(x_1,\ldots,x_{k-1},1,x_{k+1},\ldots,x_n)] . 
\end{equation}

The notion of influence is closely related to classical 
notions of voting powers.
The Banzhaf power index of $k$ in $f$ is $I_k^{\mu_{1/2}}(f)$ and 
the Shapley--Shubik power index of $k$ in $f$
is, by a theorem of Owen \cite {O},  
$\int_0^1 I_k^{\mu_p}(f)dp$. 
In \cite {GKT}  the authors proposed to define the voting power
as the probability to be pivotal 
based on realistic assumptions on individual voting distributions, and discuss
advantages and drawbacks of this approach. 

For product probability spaces, 
results of Russo, Talagrand, Friedgut and Kalai assert that
for every $p>1/2$  
sufficiently small
influences suffice to guarantee that $\mu_p(f)$ is close to 1. The latest 
such result is the following. 

\begin {theorem}[Kalai, \cite {K}] 
\label {t:cjt} 
Let $f$ be a monotone antisymmetric Boolean function. 
For every $p>1/2$ and $\epsilon$ there is  $\delta $ 
such that if $I_k^{\mu_p}(f) < \delta$ for every $k$ then 
$\mu_p(f) \ge 1-\epsilon$.
\end {theorem}

\noindent
{\bf Remark.}
The conclusion of Theorem \ref {t:cjt} remains valid if 
we replace $I_k^{\mu_p}(f)$ by the Banzhaf power index of $k$ in $f$
or by the Shapley--Shubik power index. For the Shapley--Shubik power index
a reverse implication also holds, see \cite {K}. 
We choose here a version which relies only 
on a single probability distribution $\mu_p$ and hence is more convenient 
for extensions to arbitrary probability distributions.

The purpose of this paper is to study extensions of the 
weak law of large numbers 
in the context of general probability distributions. 
Let $\mu$ be a probability distribution on $\cube$.
When $\mu$ is not a product measure
the notion of influence can be extended in different way compared to
the above.
Define the {\it effect} of the $k$'th variable on the Boolean function
$f$ as the difference between the expected value of 
$f(x_1,\dots,x_n)$ conditioned on $x_k=1$ 
and the expected value of $f(x_1,\dots,x_n)$ conditioned on $x_k=0$, and 
denote by $e_k^{\mu}(f)$ the effect of the $k$'th variable for
the Boolean function $f$, w.r.t. the distribution $\mu$. 
More precisely,
\begin{equation} \label{eq:efdef}
e_k^{\mu}(f) =  
\mu[f(X_1,\ldots,X_n) | X_k = 1] - \mu[f(X_1,\ldots,X_n) | X_k = 0].
\end{equation}
The effect is undefined if the probability for $X_k=1$ is $1$ or $0$.
Writing $\mu[X_k] = p$ and $Y_k = X_k - p$, we get
\begin{eqnarray*}
\Cov_{\mu}[f(X_1,\ldots,X_n),X_k] &=& \mu[f(X_1,\ldots,X_n) Y_k] \\ &=&
p \mu[(1-p) f | X_k = 1] + (1-p) \mu[-p f | X_k = 0] \\ &=& 
p(1-p) e_k^{\mu}(f)
\end{eqnarray*} 
so that the effect may be interpreted as a normalized form of 
the correlation between the individual vote and the election's outcome.

When $\mu$ represent a product probability measure 
(\ref{eq:product_measure}), the effect (\ref{eq:efdef}) and the 
influence (\ref{eq:Idef}) coincide, but in general this is not the case.
For instance, for general $\mu$ the effect may be negative (see item (i)
in Section \ref{s:d}) while the influence is of course always non-negative. 

It is not true that for general probability distributions
and general $f$, small influences
implies aggregation of information. Our main result is that small effects
implies aggregation of information for 
the particular case of weighted majority functions.
Moreover, unlike in Theorem \ref {t:cjt}, the bounds in our main result are 
rather realistic. 

We call monotone antisymmetric function $f$ a {\em weighted majority} function  
if there exists  non-negative weights $w_1,\dots,w_n$, not all zero such 
that $f(x_1,\dots,x_n)=1$ if $\sum_{i=1}^n w_i (2 x_i - 1) > 0$ and 
and $f(x_1,\dots,x_n)= 0$ if $\sum_{i=1}^n w_i (2 x_i - 1) < 0$.
If $n$ is odd and $w_i=1$ for every $i$, $f$ is called the majority 
function (or simple majority). 

Note that in our definition of a weighted majority function, 
if $\sum w_i (2 x_i - 1) = 0$ then 
the value of $f(x)$ may be either $0$ or $1$ as long as $f$ is monotone and 
anti-symmetric. This is different from the traditional definition of a 
weighted 
majority (or threshold) 
function where $f(x)=1$ iff $\sum w_i (2 x_i - 1) > 0$ 
and $f(x) = 0$ iff $\sum w_i (2x_i - 1) < 0$.

Thus for example, any monotone anti-symmetric 
function $f : \cube \to \{0,1\}$ 
satisfying $f(x) = 1$ when $x_1 = x_2 = 1$ and 
$f(x) = 0$ when $x_1 = x_2 = 0$ is a weighted majority function (taking
$w_1 = w_2 = 1$ and $w_3 = \cdots = w_n  = 0$) according to our definition.

The above example demonstrates that under our definition of weighted majority 
functions, there are at least $2^{2^{n-2}}$ weighted majority functions. 
Under the traditional definition the number of weighted majority 
functions is at most $2^{n^2}$ \cite{LewisCoates67,RoychSiu91}.   

\begin{figure}[htb] 
\begin{center}
\input{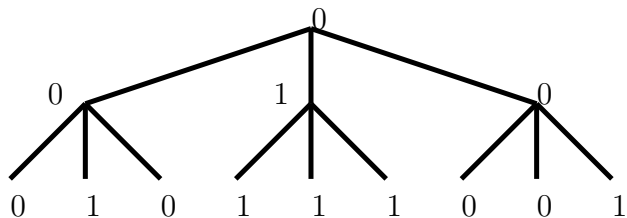} 
\caption{The function $RM_{3,2}$} 
\label{recmaj}
\end{center}
\end{figure}

Of particular interest are voting schemes 
where all the voters have the same power. 
One such case is when $f$ is invariant under a transitive group of 
permutations. In other words there exists a group of permutation 
$\Gamma \subset S_n$ such that 
$f(x_1,\ldots,x_n) = f(x_{\sigma(1)},\ldots,x_{\sigma(n)})$ for all 
$\sigma \in \Gamma$ and for all $1 \leq i,j \leq n$ there exists 
$\sigma \in \Gamma$ such that $\sigma(i) = j$; here $S_n$ denotes
the full permutation group on $n$ elements. 
One instructive example is 
the simple majority function when $n$ is odd which is invariant under $S_n$;
another is the recursive majority function $RM_{k,\ell}$
which is defined for $n=k^{\ell}$ where $k$ is odd. The definition is by 
induction. $RM_{k,1}$ is just the majority function on $k$ bits and 
\[
RM_{k,\ell+1}(x_1,\ldots,x_{k^{\ell+1}}) = 
RM_{k,1} \left(RM_{k,\ell}(x_1,\ldots,x_{k^{\ell}}),\ldots,
RM_{k,\ell}(x_{k^{\ell} - k^{\ell-1}+1},\ldots,x_{k^{\ell}}) \right).
\] 
See Figure \ref{recmaj}.

\begin {theorem}
\label {t:hkm}
$\mbox{ }$
\begin{description}
\item{\rm(a)}
For every $p > \frac{1}{2}, \epsilon >0$ 
there is $\delta = \delta (p,\epsilon) > 0$ 
such that for every weighted 
majority function $f$ and any distribution 
$\mu$ on $\cube$, if
$e_k^{\mu}[f] \leq \delta$ and $\mu[X_k = 1] \geq p$ for all $k$ 
then $\mu[f] \geq 1 - \eps$. 

In other words, if the effect of each
variable is at most $\delta$ and the probability that each variable is
$1$ is at least $p$, then $f=1$ with
$\mu$-probability at least $1-\eps$.  
\item{\rm(b)}
If $f$ is a monotone anti-symmetric function but not a 
weighted majority function, 
then there exists
a probability distribution $\mu$ such that $\mu[X_k = 1]>1/2$ 
for all $k$, yet $\mu[f] = 0$ and $e_k^{\mu}(f)=0$ for all $k$.

In other words, if $f$ is not a weighted majority function, then there
is a probability measure $\mu$ for which $f=0$ with $\mu$-probability
$1$, yet $\mu[X_k = 1] > \frac{1}{2}$ for all $k$. (Since $f$ is constant
according to the measure $\mu$, all the effects are $0$ in this case.)
\item{\rm(c)}
If $f$ is monotone anti-symmetric and invariant under a transitive group,
but is not the (simple) majority function, then then there exists 
a probability distribution $\mu$ such that $\mu[X_k = 1]>1/2$ 
for all $k$, yet $\mu[f] = 0$ and $e_k^{\mu}(f)=0$ for all $k$.
\end{description}
\end {theorem}

The rest of this paper is organized as follows.
In Section \ref {s:d} we will discuss the notions of aggregation of 
information, 
influences and effects for general probability distributions on $\cube$. 
We will try to examine what aggregation of information means 
when we do not suppose that the probability distribution for the 
voter's behavior is a product distribution. We also examine to 
what extent our technical notion of ``effects'' represent real 
influence in the non-technical sense of the words. 
Section \ref {s:p} contains the proof of our theorem  
and in Section \ref {s:prob} we present 
several natural problems as well as an 
example showing that 
Theorem \ref {t:cjt} does not extend to arbitrary Boolean 
monotone functions even for the restricted class of FKG-distributions.
Finally, in Section \ref{sect:appendix}, we present an alternative
proof of Theorem \ref{t:hkm} (a) that yields sharper quantitative bounds.

\section {Voting games, information aggregation and notions of influence}
\label {s:d}

Consider 
the following scenario. Every agent $k$ receives a single bit of
information $s_i$ which is either `Vote for Alice' or `Vote for Bob' and
these signals are independent. When Alice is the better candidate 
the probability of receiving  the signal `Vote for Alice' is
$p>1/2$. 
Condorcet's Jury Theorem deals
with the case that the voters vote precisely as the signal dictates and the
decision is made according to  the simple majority rule. It asserts that
for every $p>1/2$  the better candidate will be elected with 
probability tending to one. Thus the majority rule allows to reveal 
the actual state of the world from rather weak individual signals.

A major problem in the economic and political
interpretation of Condorcet's Jury Theorem and its extensions
arises from the fact that
the basic assumption of probability
independence among voters
is quite unrealistic.
Without the assumption of independence, Condorcet's Jury Theorem as stated
is no longer true, and it will no longer be the case that
when each individual votes for Alice with probability $p>1/2$,
Alice will win with a high probability.

To see this, consider the following
example. As before, we have an election between Alice and Bob and
Alice is the superior candidate. The distribution of signals
$s_1,s_2, \dots , s_n$ will be biased towards Alice as follows:
Let $p=1/2+\epsilon/2$, where $\epsilon$ is small.
First choose at random a number $t$ uniformly
in the interval $[\epsilon,1]$.
Then, independently for each $i$, 
choose  the $i$'th voter  signal $s_i$ to be `1' with
probability $t$ and `0' with probability $1-t$. Voters with $s_i=1$
will vote for Alice.
In this case, the probability
for each individual signal $s_i$ being `1' is $p$
but the individual signals are not independent. The probability that
Alice will win is below $\frac{1}{2(1-\eps)}$  for any
number of voters. This is because we can think of $t$ being chosen in two 
stages. First we toss a coin which is 'H' with probability $\epsilon/(1-\epsilon)$. If the 
coin is 'H', then $t$ is chosen uniformly in the interval $[1-\epsilon,1]$. 
This contributes to the probability that Alice wins at most 
$\epsilon/(1-\epsilon)$. If the coin is 'T' then 
$t$ is chosen uniformly in the 
interval $[\epsilon,1-\epsilon]$. Here by symmetry, Alice and Bob have the 
same chance of winning. Thus the contribution to the probability that Alice 
will win from this case is $\frac{1-2\epsilon}{2(1-\epsilon)}$. Thus the 
overall probability that Alice will win is at most 
$\frac{1-2\epsilon}{2(1-\epsilon)} + \frac{\epsilon}{1-\epsilon} = 
\frac{1}{2(1-\eps)}$.

An even more extreme example is the case in which all voters vote in the same 
way: With probability $p$ they all vote for Alice and with probability $1-p$ 
they all vote for Bob. Alice
will be elected with probability $p$ regardless of the number of
voters when the election is based on simple majority and for every
other simple game.

These simple examples will help us to
examine the notions of information aggregation and influence 
in the case when the assumption of probability independence is dropped.
The problem in these examples is not
in the way information aggregates but in
the quality of the information to start with.
This assertion can be formalized as follows.
Suppose that Alice and Bob are given an a-priori probability $1/2$ of being 
the superior candidate. We assume that the distribution of voters for Bob 
given Bob is the superior candidate is the same as the distribution of 
voters for Alice given that Alice is the superior candidate. Thus in the 
first example above if Bob is superior then we choose $t$ uniformly in 
$[\eps,1]$ and then each voter votes for Bob independently with 
probability $t$. In the second example if Bob is superior, then all voters
will vote for Bob with probability $p$.
 
We now wish to decide between the hypothesis that Alice is the superior 
and the hypothesis that Bob is the superior candidate  
given the entire vector of individual signals. 
It is intuitively clear and easy to prove using the Neyman-Pearson Lemma 
that in both cases described above 
one should guess that Alice is superior to Bob exactly when the majority 
of voters voted for Alice. However, in both examples above the probability 
that the majority will vote Alice when Alice is superior 
is bounded away from $1$ and
tends to $1/2$ as $p$ does.

When we consider general distributions, the issue is to
understand what information we can derive on the superior alternative
from knowing the signals of all individuals and how the voting mechanism
extracts this information. Note that in the examples we considered 
above the individual effects are large while the individual 
influences are small. This is most transparent in the second example where 
if $f$ is the majority function and $n \geq 3$, 
then all of the influences are $0$, while the effect of 
all voters are $1$. 
Theorem \ref {t:hkm} asserts that for the weighted 
majority voting rule (and only for these rules) 
for every probability measure on $\cube$,
small individual effects 
implies asymptotically complete aggregation of information. 

Let us next consider the notion of influence without probability independence.
The notion of pivotal variables (or players) and 
influence is of important technical importance in various 
areas of mathematics, computer science and economics. 
This notion is also of a considerable conceptual importance. 
The voting power index of Banzhaf is 
based on measuring the influence with respect to the uniform 
distribution. The Shapley--Shubik power index can also be based on 
the influence with respect to another distribution.
Conceptual understanding of voting power in situations
where the voters' behavior is not independent is of great interest.
In \cite {GKT}  the authors propose to define the voting power
as the probability to be pivotal 
based on realistic assumptions on individual voting distributions.
We make the following remarks on the notion of individual effects which is 
quite a different extension of influence and voting power measures
to arbitrary probability distributions. 




\begin{description}
\item{(i)} 
For general distributions, the effect of an agent can be negative. 
This will be the case for a voter who always votes  for
the candidate who is the  underdog in the election polls and also for a 
committee member who antagonizes the other members of the committee.
(On the other hand, the influence of an agent is always nonegative,
because it is defined as a represents a probability.) 
\item{(ii)} 
A dummy (a voter $k$ which is never pivotal)
has zero influence (with respect to every probability distribution). 
He may nevertheless have a large effect, such as if 
he always votes  for the candidate who is expected
to win  according to election polls. 
In real life, this will also be the case
for an observer on a committee
without the right to vote but who is likely to convince the
committee of  his opinion. Note that in
the first case we do  not attribute to that player real
``influence'' in the (non-technical) English sense of the word,
while  in the second case we would consider him ``influential''.
The uncertainty in  interpreting  effects as real
``influences'' is genuine. 
\item{(iii)} 
What is the 
motivation for a voter to vote, given the small probability for 
him to be pivotal? This is a social dilemma, related to, e.g., 
the so-called {\em tragedy of the commons}, and has been extensively 
discussed in the political science and philosophy literature. (Sometimes
the term {\em voting paradox} has been used for this dilemma, but
may cause some confusion as
the same term is used also for Condorcet's famous observation that
when three or more choices are available, the majority preference
between them need not be transitive.)
A possible solution to the dilemma
may lie in the fact that in real-life elections, individual effects 
tend to be large, namely bounded away from zero regardless of the 
size of the society. 
The uncertainty in 
regarding effects as real ``influence'' may suggest that it is 
the effect of an agent rather than his influence which is related
to his ``satisfaction'' with  the social decision process and his
ability to identify with the collective choice. 
\end{description}



\section{Proof of Theorem \ref {t:hkm}}
\label {s:p}

\subsection{Part (a) of the theorem}

We begin this section by providing a probabilistic proof of the
following result, which clearly implies 
Theorem \ref{t:hkm} (a). 

\begin{lemma} \label{lem:prob}
Let $(w_i)_{i=1}^n$ be non-negative weights which are not all $0$, 
let $0 < q < 1$, and let $f : \cube \to \bits$ be a function which satisfies 
\[
f = \left\{
\begin{array}{ll}
1 & \mbox{if } \sum_{i=1}^n (2 x_i - 2 q) w_i > 0 \\
0 & \mbox{if } \sum_{i=1}^n (2 x_i - 2 q) w_i < 0.
\end{array}  \right. 
\]
Write $W = \sum_{i=1}^n w_i$.
Suppose furthermore that $p > q$ and that $\mu$ is a 
probability measure satisfying 
$\mu[X_i] = p_i$ and
\begin{equation}  \label{eq:first_cond}
\sum_{i=1}^n w_i p_i = p W
\end{equation}
as well as
\begin{equation}  \label{eq:second_cond}
\sum_{i=1}^n w_i p_i (1-p_i) e_i^{\mu}[f] \leq p (1-p) \delta W . 
\end{equation}
Then 
\[
\mu[f] \geq 1 - \frac{\delta p(1-p)}{p-q}.
\]
\end{lemma}
(Note that (\ref{eq:first_cond}) holds if $\mu[X_i] = p$ for all $i$,
and that (\ref{eq:second_cond}) holds if 
$\mu[X_i] = p$ and $e_i^{\mu}[f] \leq \delta$ for all
$i$, so that indeed Theorem \ref{t:hkm} (a) follows.)

\begin{proof}
Let $X = \sum_{i=1}^n (2 X_i - 2 q) w_i$.
We start by noting that $\mu[X] = (2p-2q)W$.

We let $g = 1 - f$ and $Y_i = p_i - X_i$, so
that 
\[
\mu[Y_i g] = \Cov_{\mu}[f,X_i] = p_i(1-p_i) e_k^{\mu}[f]. 
\]

Note that conditioned on $g = 1$, 
$\sum_{i=1}^n (2 X_i - 2 q) w_i \leq 0$ and therefore 
$\sum_{i=1}^n w_i Y_i \geq  (p-q) W$. 
It follows that 
\begin{eqnarray}
\nonumber
\mu\left[\left(\sum_{i=1}^n w_i Y_i\right) g(X_1,\ldots,X_n)\right] 
& \geq & (p-q) W \mu[g] \\
& = & (p-q) W (1 - \mu[f]).   \label{eq:g1}
\end{eqnarray}
On the other hand,
\begin{eqnarray}
\nonumber
\mu\left[\left(\sum_{i=1}^n w_i Y_i\right) g(X_1,\ldots,X_n)\right] 
& = & \sum_{i=1}^n w_i \mu[Y_i g(X_1,\ldots,X_n)] \\
\nonumber
& = & \sum_{i=1}^n w_i p_i (1-p_i) e_i^{\mu}[f] \\
& \leq  & p(1-p) \delta W. 
 \label{eq:g2}
\end{eqnarray}
Combining (\ref{eq:g1}) and (\ref{eq:g2}), we get that  
\begin{eqnarray*}
\mu[f] & \geq & 1 - \frac{\sum_{i=1}^n w_i p_i (1-p_i) e_i^{\mu}[f]}{(p-q) W} 
\\ & \geq &  1 - \frac{\delta p(1-p)}{p-q}.
\end{eqnarray*}
\end{proof}

\subsection{Parts (b) and (c) of the theorem}

We note that part (c) of Theorem \ref{t:hkm} follows immediately from part
(b), because the only weighted majority function that is
invariant under a transitive group, is simple majority. 
Let us nevertheless begin by giving an 
independent and simple proof of part (c).  
Note that if $f$ is not the majority function 
then there is a vector $(x_1,x_2,\dots,x_n) \in \{0,1\}^n$ such that 
$f(x)=0$ and $x_1+x_2+\dots+x_n>n/2$. Then we can simply take $\mu$ to be 
uniform probability distribution on the orbit of $x$ under $\Gamma$.
It is then easy to see that $\mu[X_k] > 1/2$ for all $k$ and that 
$\mu[f=0]=1$. 

We now turn to the proof of Theorem \ref{t:hkm} (b).
We will show that if $f$ is not a weighted majority function, then 
there exists a measure $\mu$ satisfying $\mu[X_k] > 1/2$ for all $k$ 
and $\mu[f=0]=1$.  


Define $[n]=\{1,2,\dots,n\}$.
For $S \subset [n]$ put $x_S=(x_1,x_2\dots,x_n)$ 
where $x_i=1$ if and only if $i \in S$. 
Let $H$ be a hypergraph whose set of vertices is 
$[n]$ and whose edges are subsets $S$ of $[n]$ such that $f(x_S)=0$. 
Let $\tau^* = \tau^*(H)$ be the fractional cover number of $H$, i.e., 
the infimum over all $\nu:\cube \to \R$ of 
$\sum_{S \in H} \nu[x_S]$, under the condition that 
$\nu(x_S) \ge 0$ for every $S \in H$ and 
$\sum_{S \in H, k \in S} \nu[x_S] \ge 1$ for all $k$. 
We get $\tau^{\ast} = \infty$ 
if there are no $\nu$ satisfying the two conditions above (note that this is 
the case if $f(x) = x_1$, say).


If $\tau^* < 2$, then we can define $\mu(S)=0$ if $f(S)=1$ and
$\mu(S)=\nu(S)/\tau^*$ when $f(S)=0$. The probability measure $\mu$ satisfies
that 
\[
\sum_{S:k \in S, f(S)=0} \mu (x_S)\geq 1/\tau^* > 1/2
\]
for every $k$ and $\mu[f=0]=1$ as stated in the theorem. 
Therefore, in order to prove part (b) of the theorem, 
it only remains to analyze the case $\tau^* \geq 2$. 
  
A well known equivalent (by linear programming duality) 
definition of of $\tau^*$ is as the supremum of $\sum_{i=1}^n w_i$ under 
the condition that $w_k \ge 0$ for $k=1,2,\dots,n$ and 
$\sum \{w_i: i \in S \} \le 1$ for every $S \in H$. 

Assume first that $\tau^{\ast} > 2$. In this case we can find $w_i$'s such 
that $\sum_i w_i > 2$ and $f(x_1,\dots,x_n) = 1$ if $\sum w_i x_i > 1$. 
By slightly perturbing the $w_i$ we may assume that for all $x \in \cube$ 
it holds that $\sum_i w_i x_i \neq \frac{1}{2}\sum_i w_i$ in addition to the 
properties that $\sum_i w_i > 2$ and $f(x_1,\dots,x_n) = 1$ if 
$\sum w_i x_i > 1$. 
Let $g(x) = 1$ if $\sum_i w_i x_i > \frac{1}{2}\sum_i w_i$ and 
    $g(x) = 0$ if $\sum_i w_i x_i < \frac{1}{2}\sum_i w_i$.
Then $g$ is 
anti-symmetric and $f=0 \implies g = 0$. It follows that $f=g$ so that 
$f$ is a weighted majority function as needed. 

The remaining case is where $\tau^{\ast} = \sum w_i=2$. 
We obtain that $f(x_1,\dots,x_n) = 1$ if $\sum w_i x_i > 1$.
Since $f$ is anti-symmetric it follows that  
$f(x_1,\dots,x_n) = 0$ if $\sum w_i x_i < 1$. 
The result follows. 
{\hfill $\square$}

\section {Problems and an additional example}
\label {s:prob}

The following problems naturally suggest themselves at this point:

\begin{description}
\item{(1)} 
For which class of distributions is it the case that for simple majority 
small voting power implies asymptotically complete aggregation of information?
\item{(2)} 
For which class of distributions is it the case that for every monotone 
Boolean function small voting power implies 
asymptotically complete aggregation of information?
\item{(3)} 
For which class of distributions is it the case that for every monotone 
Boolean function small individual effects 
implies asymptotically complete aggregation of information?
\end{description}


A natural condition to impose on the 
distribution $\mu$ which is realistic in various economic 
situations is the FKG
condition (see \cite {L}).
For $x=(x_1,\dots,x_n)$ and $y=(y_1,\dots,y_n)$, define 
\[
\max (x,y)=(\max (x_1,y_1),\dots, \max (x_n,y_n))
\]
and
\[
\min (x,y)=(\min (x_1,y_1),\dots,\min (x_n,y_n)) . 
\]
One definition of FKG measure on $\{0,1\}^n$
goes as follows:
A distribution $\mu$ on $\{0,1\}^n$ (or on $\R^n$)
is called an FKG measure if
for every $x,y \in \{0,1\}^n$ we have
\[
\mu(x)\mu(y) \leq \mu(\max (x,y)) \mu (\min (x,y)).
\]
The FKG property is a profound notion
of non-negative correlations between agents' signals. It implies
(but is strictly stronger than) the following condition
(known as {\it non-negative association}, see \cite {M}):
For all increasing
real functions $f$ and $g$, it is the case that $E[fg] \geq E[f] E[g]$.
This is equivalent to the condition that for all increasing events $A$ and
$B$ we have that $P[AB] \ge P[A]P[B]$. Under the FKG property
if the simple game is monotone,  all effects are
non-negative. This form of non-negative correlation is a plausible
assumption to make in various contexts of collective choice.
It is easy to see that under the condition of 
non-negative association all individual effects are non-negative.

\begin{description}
\item{(4)} 
For which class of monotone Boolean functions does small individual effects 
imply asymptotically complete aggregation of information?
\end{description}

In the following subsection, 
we present an example of an FKG measure and a 
monotone Boolean function such that the 
individual effects are small and yet there is no asymptotically 
complete aggregation of information. In this example both the voting scheme 
and the measure $\mu$ are invariant under a transitive group of permutations. 

\subsection{Example: FKG without aggregation}
{\bf The measure $\mu$.}
We start by describing the measure $\mu$. The measure is given 
by a Gibbs measure for the Ising model on the $3$-regular tree.
See e.g. \cite{G,P}.  
The measure is defined as follows. Let $T_r = (V_r,E_r)$ 
be the $r$-level $3$-regular tree.  
This is a rooted tree where each internal nodes has exactly $3$ children 
and all the leaves are at distance exactly $r$ from the root $\rho$. 
Let $L_r$ be the set of leaves of that tree. Note that $|L_r| = 3^r$.  
Thus in Figure \ref{recmaj} the underlying tree is $T_2$. 

We first define a measure $\nu$ on the tree $\{0,1\}^{V_r}$. 
In this measure the probability of $x$ is given by 
\[
\nu[x] = \frac{1}{2} 
\prod_{(u,v) \in E_r} \left(
(1-\epsilon) 1_{\{x_u = x_v\}} + \epsilon 1_{\{x_u \neq x_v\}}
\right).
\]
In words, this means
that in the measure $\nu$ the sign of the root $x_{\rho}$ is chosen 
to be $0$ or $1$ with probability $1/2$. Then each vertex
inherits its parents label with probability $\theta=1-2\eps$ and is
chosen independently otherwise. 

Our measure $\mu$ is defined on $\{0,1\}^{L_r}$ (so that
the voters are the leaves of the tree) 
as follows. 
\[
\mu[x] = \sum_{y : y | L_r \leq x} 
\nu[y] \delta^{|\{i : x_i = 1, y_i = 0\}|}.
\]   
In other words, a configuration of votes according to $\mu$ may be obtained 
by drawing a configuration $x$ according to $\nu$ and looking at $x | L_r$. 
Then for each of the coordinates of $i \in L_r$ independently, the vote at 
$x$ re-sampled to have the value $1$ with probability $\delta$. Below we will 
sometime abuse notation and write $\mu$ for the joint probability 
distribution of $x$ and $y$. 

Standard results for the Ising model (see, e.g., \cite{G}) imply 
that $\mu$ is an FKG measure. 
Moreover it is easy to see that the measure is 
invariant under a transitive group and  
that $\mu[x_i] = (1+\delta)/2$ for all $i$.

{\bf The function $m$.} 
The function $m$ is given by the recursive majority function 
$m=RM_{3,r}$. Clearly, $m$ is monotone, anti-symmetric and invariant under 
a transitive group. 

\begin{claim}
If $\eps = \delta \leq 0.01$ then  
$\mu[m] \leq 1/2 + \delta/2$ for $m=RM_{3,r}$ and all $r$. 
\end{claim}  

\begin{proof}
The proof below is similar to arguments in \cite{M1,M2}. 
Let $(y_v : v \in V_r)$ 
be chosen according to the measure $\nu$. Let $(x_v : v \in L_r)$ be 
obtained from $y_v$ by re-sampling each of the coordinates of 
$(y_v : v \in L_r)$ to $1$ with probability $\delta$. Let $(m_v : v \in V_r)$
denote the value of the recursive majority of all $(x_w : w \in L_r(v))$, 
where $L_r(v)$ are all the leaves of $T$ below $v$.
We will show that $\mu[m = m_{\rho} = 0 |  y_{\rho}=0] \geq 1-\delta$. 
Since $\mu[y_{\rho}=0] = 1/2$,  
we conclude that 
$\mu[m] \leq 1/2 + \mu[m | x_{\rho} = 0]/2 \leq (1+\delta)/2$, as needed. 

We are interested in the probability that $m_v = 0$ conditioned 
on $y_v = 0$. It is easy to 
see that this probability only depends on the height of $v$, i.e., on the
distance between $v$ and the set of leaves. We let $p(k)$ denote the 
probability that $m_v = 0$ conditioned on $y_v = 0$ 
for a vertex $v$ of height $k$. 

Clearly, $p(0) = 1-\delta$. We want to prove by induction that 
$p(k) \geq 1-\delta$ for all $k$. Let $v$ be a node of height $k+1$ and $w$ a 
child of $v$. Note that conditioned on $x_v = 0$ the probability that 
$m_w = 0$ is at least $(1-\eps)p(k)$ which is at least 
$t = (1-\eps)(1-\delta)$ by the induction hypothesis. 
Moreover, noting        
that the values of the
majorities of the children of the node $v$ are conditionally independent
given that $m_v = 0$, we conclude that 
\[
p(k) \geq t^3 + 3t^2(1-t) = 3t^2 - 2t^3 = t^2(3-2t).
\]
We need that $t^2(3 - 2t) \geq 1-\delta$ or recalling that $\eps=\delta$: 
$(1-\eps)^4 (3 - 2(1-\eps)^2) \geq (1-\eps)$. This in turn is equivalent to 
$(1-\eps)^3(3-2(1-\eps)^2) \geq 1$. The function $h(\eps) = 
(1-\eps)^3(3-2(1-\eps)^2)$ has $h'(\eps) = 10(1-\eps)^4 - 9(1-\eps)^2 = 
(1-\eps)^2(10(1-\eps)^2 - 9)$. Therefore $h$ is increasing in the interval 
$[0,0.01]$. Since $h(0)=1$ it follows that $h(\eps) \geq 1$ for all 
$\eps \leq 0.01$ as needed.
\end{proof}

Our next objective is to bound the effect of a voter at level $r$. 
We will prove the following: 

\begin{claim} \label{cl:1}
The measure $\mu$ on $T_r$ and the function $m=RM_{3,r}$ 
satisfy that the effect of each voter is at most 
$(1-\eps/2)^{(r-1)/2} + 2^{-(r-1)/2}$.
\end{claim}

\begin{proof}
The argument here is similar in spirit to an argument in \cite{BKMP}. 
Let $t+s=r$ where $t \geq (r-1)/2$ and $s \geq (r-1)/2$. 
Fix a leaf voter $i$. We want to estimate 
$\mu[m = 1 | y_i = 1] - \mu[m = 1 | y_i = 0]$. Let's denote 
by $\mu_0$ the measure $\mu$ conditioned on $y_i=0$ and by $\mu_1$ the
measure $\mu$ conditioned on $y_i=1$. 

Let $i=v_0,v_1,\ldots,v_r = \rho$ denote the path from $i$ to the root. 
We first claim that the measures $\mu_0,\mu_1$ and $\mu$  
may be coupled in such a way 
that except with probability $(1-2\eps)^t$ the only disagreements between 
$\mu_0, \mu_1$ and $\mu$ are on vertices below $v_t$.

The follows immediately from the random cluster representation of the model.
In this representation we declare and edge $(u,v)$ open with probability 
$(1-2\eps)$ and closed with probability $2\eps$. If the edge $(u,v)$ 
is open then $y_u = y_v$, otherwise, 
the two labels are independent. It is then clear that 
we may couple the two measure $\mu_0,\mu_1$ and $\mu$ 
below $v_t$ as long as the 
path from $i$ to $v_t$ contains at least one closed edge. The probability that 
such an edge does not exist is at most $(1-2\eps)^t$. The proof of the first 
claim follows. 

For each $j$ denote by $u_j$ and $w_j$ the siblings of $v_j$. We assume  
that the measures $\mu_0,\mu_1$ and $\mu$ are coupled in such a way 
that the only disagreements between them are on vertices below $v_t$. 
Note that if this is the case, then if the values of $m$ under 
$\mu_0$ and $\mu_1$ are different 
then for all $r \geq j \geq t$ it holds that $m_{u_j} \neq m_{w_j}$. 
We wish to bound the $\mu$ probability that $m_{u_j} \neq m_{w_j}$ 
for $r \geq j \geq t$. We will bounds this probability conditioned 
on the values $(y_{v_j})_{j=t}^r$. Conditioned on $(y_{v_j})_{j=t}^r$ 
the event $m_{u_j} \neq m_{w_j}$ are independent for different $j$'s. 
Moreover, by the Markov property, 
$\mu[m_{u_j} \neq m_{w_j} | (y_{v_h})_{h=t}^r] = 
\mu[m_{u_j} \neq m_{w_j} | y_{v_{j-1}}]$. Finally note that conditioned 
on $y_{v_{j-1}}$ the random variables $m_{u_j},m_{w_j}$ are identically 
distributed and independent. Therefore 
\[
\mu[m_{u_j} \neq m_{w_j} | y_{v_{j-1}}] \leq 
\max_{p \in [0,1]} 2p(1-p) \leq 1/2.
\]
We thus obtain that the $\mu$ probability that $m_{u_j} \neq m_{w_j}$ 
for $r \geq j \geq t$ is at most $2^{-s}$.
\end{proof}


%


\section{A linear programming bound}
\label{sect:appendix}

We present in this section an alternative proof of Theorem \ref{t:hkm}
(a) using linear programing. This approach
yields tight bounds, stated in the following lemma. 
\begin{lemma} \label{lem:lin}
Let $(w_i)_{i=1}^n$ be positive weights, let $0 < q < 1$, and let 
$f : \cube \to \bits$ be a function satisfying 
\[
f = \left\{
\begin{array}{ll}
1 & \mbox{if } \sum_{i=1}^n (2 X_i - 2 q) w_i > 0 \\
0 & \mbox{if } \sum_{i=1}^n (2 X_i - 2 q) w_i < 0 .
\end{array} \right. 
\]
Write $W = \sum_{i=1}^n w_i$.
Suppose that $p > q$ and that $\mu$ is a probability measure satisfying
\begin{equation}   \label{eq:3rd_cond}
\sum_{i=1}^n \mu[X_i] = p W
\end{equation}
and 
\begin{equation}  \label{eq:4th_cond}
\sum_{i=1}^n w_i p_i (1 - p_i) e_i^{\mu}[f] \leq \delta W p (1-p).
\end{equation}
If $\delta \geq \frac{p-q}{p(1-q)}$, then 
\[
\mu[f] \geq \frac{p-q}{1-q}, 
\]
whereas otherwise
\[
\mu[f] \geq \max\left\{\delta p, 1 - \frac{\delta p (1-p)}{p-q}\right\}.
\]
These bounds are tight.
\end{lemma}
(Note that the conditions of Theorem \ref{t:hkm}
(a) imply (\ref{eq:3rd_cond}) and (\ref{eq:4th_cond}).)

\begin{proof}
We will first make the necessary computations for the case of simple majority.
Let $f : \{0,1\}^n \to [0,1]$ be a symmetric monotone function.
Let $\delta_i(x) = 1$ if $x_i = 1$ and $\delta_i(x) = 0$ otherwise.
Let $\eta_i(x) = 1 - p$ if $x_i = 1$ and $\eta_i(x) = -p$ if $x_i = 0$.

We want to minimize 
\begin{equation} \label{eq:cost}
\mu[f(X_1,\ldots,X_n)] = \sum_{x} \mu(x) f(x).
\end{equation}
under the restrictions that
\begin{equation} \label{eq:const1}
\sum_{i,x} \delta_i(x) \mu(x) = \sum_{i=1}^n \mu[X_i] = n p,
\end{equation}
and (letting $Y_i = X_i - p$)
\[
\sum_{i=1}^n \mu[Y_i f(X_1,\ldots,X_n)] = 
\sum_{i=1}^n \Cov_{\mu}[f(X_1,\ldots,X_n),X_i] \leq \delta p(1-p) n,
\]
which gives
\begin{equation} \label{eq:const2}
\sum_{x,i=1}^n \eta_i(x) \mu(x) f(x) \leq n \delta p (1 - p).
\end{equation}

The constraints (\ref{eq:const1}) and (\ref{eq:const2}) and the
cost (\ref{eq:cost}) are invariant under the action of $S_n$ on the
coordinates of $x$ (since $f(x)$ has this
invariance property). It follows that there exists a minimizer $\mu$
which is symmetric, i.e.,
\[
\mu(x) = \frac{a_{|x|}}{\binom{n}{|x|}},
\]
where $a$ is a positive function.

Since $f$ is a majority function, it follows that there exists an $r$ 
such that $f(x) = f(|x|) = 1$ if and only if $|x| > r$. We
therefore obtain the following optimization problem. Write $q = r/n$
and $q' = (r+1)/n$. We assume below that $p > q'$.

We want to minimize
\begin{equation} \label{eq:n_cost}
\sum_{i=r+1}^n a_i,
\end{equation}
under the restrictions
\begin{equation} \label{eq:n_const1}
\sum_{i=0}^n a_i \frac{i}{n} = p,
\end{equation}
and
\begin{equation} \label{eq:n_const2}
\sum_{i=r+1}^n a_i \frac{i}{n} - p \sum_{i=r+1}^n a_i 
\leq \delta p (1-p).
\end{equation}

It is useful to introduce $A = \sum_{i=r+1}^n a_i$ and $B =
\sum_{i=r+1}^n a_i \frac{i}{n}$. Similarly, we write $A' =
\sum_{i=0}^r a_i$ and $B' = \sum_{i=0}^r a_i \frac{i}{n}$. 
Note that $A,A',B,B'$ are all positive, $0 \leq B' \leq qA'$,
similarly, $q'A \leq B \leq A$. Note that $A + A' = 1$, 
(\ref{eq:const1}) may be written
as $B + B' = p$ and (\ref{eq:const2}) may be written as 
$B - p A \leq \delta p (1-p)$. 

Moreover, it is easy to see that any $A,A',B,B'$ which satisfy the
above equations give rise to $a_i$ satisfying the constraints. 
Using $A+A'=1$ and $B+B'=p$ we thus led to the following optimization
problem in $A$ and $B$: Minimize $A$ under the constraints
\begin{equation} \label{eq:AB1}
0 \leq A \leq 1,
\end{equation}
\begin{equation} \label{eq:AB2}
q'A \leq B \leq A,
\end{equation}
\begin{equation} \label{eq:AB3}
0 \leq (p - B) \leq q (1 - A),
\end{equation}
\begin{equation} \label{eq:AB4}
B - p A \leq \delta p (1-p).
\end{equation}

In other words, we are looking for the minimal $A$ satisfying
\[
\begin{array}{cc}
\max\{B,\frac{B}{p} - \delta(1-p)\} \leq A \leq 
\min\{\frac{B}{q'},1-\frac{p}{q}+\frac{B}{q}\}, & 0 \leq B \leq p.
\end{array}
\]
From the assumption $p > q'$ it follows that
for $0 \leq B \leq p$, the minimum on the right hand side is obtained
by $1 - \frac{p}{q} + \frac{B}{q}$. We may thus simplify:
\[
\begin{array}{cc}
\max\{B,\frac{B}{p} - \delta(1-p)\} \leq A \leq 
1-\frac{p}{q}+\frac{B}{q}, & 0 \leq B \leq p.
\end{array}
\]
The two functions bounding $A$ from below are increasing in $B$. 
Therefore in order to minimize $A$ we should minimize $B$. 
Note that $B > \frac{B}{p} - \delta(1-p)$ if and only if 
$B < B_c = \delta p$. 

Suppose that $B < B_c$. Then we obtain that 
$B \leq 1 - \frac{p}{q} + \frac{B}{q}$, or equivalently, 
$B \geq \frac{p-q}{1-q}$. We thus conclude that if $\frac{p-q}{1-q}
\leq \delta p$, then the minimum for $B$ (and therefore for $A$) is
obtained at $\frac{p-q}{1-q}$. 

If $\frac{p-q}{1-q} > \delta p$, then $B > B_c$. Moreover, in this
case, $\frac{B}{p} - \delta(1-p) \leq 1-\frac{p}{q}+\frac{B}{q}$, 
which is equivalent to 
\[
B \geq p - \frac{\delta p q (1-p)}{p-q}.
\]
Therefore the minimal $B$ in this case is given by 
\[
B = \max\left\{\delta p,p - \frac{\delta p q (1-p)}{p-q}\right\}
\]
and therefore by the bound $A \geq \frac{B}{p} - \delta(1-p)$ 
the minimal $A$ is given by 
\[
A = \max\left\{\delta p, 1 - \frac{\delta p (1-p)}{p-q}\right\} \, . 
\]
This establishes the lemma for the special case of simple majority. 
Moving from simple majority to weighted majority is easy. 
First note that we
can assume that all weights are nonnegative integers. Replace a  variable 
$x_k$ with $w_k$ copies. We will thus consider $\{0,1\}^W$ where 
$W=w_1+w_2+\dots w_n$. Consider the distribution $\nu'$ on $\{0,1\}^W$ 
induced from $\nu$ with the requirement that with probability 1
all copies of an original variable have the same value.
The desired result for the weighted majority on $\{0,1\}^n$ follows from 
the case of simple majority on $\{0,1\}^W$. 
\end{proof}

\bigskip\noindent
{\bf Acknowledgement.} The 
authors are grateful to Yuval Peres for fruitful discussions. 


\end{document}